\documentstyle [12pt,twoside, epsf]{article}

\def\picill#1by#2(#3)
{\vbox to #2
{\hrule width #1 height 0pt depth 0pt
\vfill\epsffile{#3}}}

\begin{document}

\date{}

\title{\bf Lune-Free Knot Graphs}

\author{
Shalom Eliahou\\
Universit\'e du Littoral C\^ote d'Opale\\
Calais, France\\
and\\
Frank Harary\\
New Mexico State University
Las Cruces, New Mexico\\
and\\
Louis H. Kauffman \\
  Department of Mathematics, Statistics and Computer Science \\
  University of Illinois at Chicago \\
  851 South Morgan Street\\
  Chicago, IL, 60607-7045
}

 \maketitle

 \thispagestyle{empty}

 \subsection*{\centering Abstract and Dedication}

{\em
 This paper is an exploration of simple four-regular graphs in the plane (i.e. loop-free and with no more than one edge
between any two nodes). Such graphs are fundamental to the theory of knots and links in three dimensional space, and their
planar diagrams. We dedicate this paper to Frank Harary (1921 -- 2005), whose fascination with graphs of knots inspired this work, and with whom we
had the pleasure of developing this paper.}

\section{Introduction}
This paper is an exploration of four-regular graphs in the plane that have no more than one edge between any two nodes.
We call such graphs ``lune-free" since the appearance of two edges from one node to another makes a lune-shape in the plane.
Such graphs are fundamental to the theory of knots and links in three dimensional space, and their planar diagrams.
\bigbreak

The reader interested in the relationship of these graphs with knot theory should consult the fundamental paper of John H. Conway
\cite{Conway}. In this paper, Conway uses a small collection of lune-free graphs to construct and enumerate all knots (knot diagrams up to 
topological equivalence) with no more than ten crossings. There is a potential for using the full class of lune-free diagrams to enable deeper
tabulations of knots and links. See \cite{Jablan} for recent work along these lines.
\bigbreak

We became interested in these graphs for their own sake, motivated by their possible use in the theory of knots and links.
In particular, seeing that the first lune-free link projection is the familiar graph of the Venn diagram for three circles, and that the first
single component (in the sense of knot theory) lune-free graph had eight nodes, we wanted to prove that beyond eight nodes, there were lune-free
graphs of one component for every choice of the number of nodes. This conjecture turned out to be true. It is a consequence of some interesting 
inductive constructions that we produce in the third section of the paper. These constructions depend upon an admissibility 
criterion described below. It was then natural to ask about the existence of inadmissible lune-free graphs, and the remainder of the paper 
is devoted to theorems describing when they can be constructed. These results are new, and they are answers to questions that arise naturally 
in this subject.
\bigbreak

Here is an outline of the contents of the paper.
We first show that, for projected knot diagrams, there are no lune-free diagrams with less than $8$ nodes, and that from then on 
there is at least one such diagram with $v$ nodes for every $v \ge 8.$ This result is proved by using the Euler formula for plane graphs and
some graphical recursions in Sections 2 and 3. A lune-free knot graph is said to be \textit{admissible} if it has two adjacent regions with 
both regions having at least 4 sides. The point about an admissible knot 
graph is that  it admits a certain recursive procedure, described in 
Section 3, that allows us to make infinitely many lune-free knot graphs from it. In the remaining sections of
the paper we explore the existence question of {\em tight}, i.e.
inadmissible, lune-free knot and link graphs.
The basic question, that we shall completely
answer here, is the following: for which integers $v$ does there
exist a {\em tight} lune-free knot or link graph, respectively?
For links, the answer is: all integers $v \geq 8$ such that
$v \not \equiv \pm 1 \bmod 6$ can be realized as the number of crossings
of such graphs. The answer is the same for knots as for links, with the sole 
exception of $v=12$, which cannot be realized by a tight lune-free {\em knot} 
graph. (However, there is a tight lune-free link graph with 12 crossings and 
$\mu = 3$ components.) This is accomplished in Section 4. Finally, in Section 5,
we construct knot graphs with a fixed number of lunes and with any sufficiently 
large number of vertices.
\bigbreak

\noindent {\bf Acknowledgement.} 
It gives the third author great pleasure to acknowledge support from NSF Grant DMS-0245588.
\bigbreak

\section {Applying the Euler Formula}

We shall sometimes use the knot theory terminology {\it universe} for a $4$-regular
connected plane multi-graph. This paper is concerned with those universes that are 
simple graphs in the sense that they have no more than one edge between two nodes,
and have no loops (a loop is an edge whose endpoints are identical). We shall call
a simple and connected universe a {\it lune-free graph}. There are other
terminologies for this concept. Such graphs are called {\it Conway polyhedra} \cite{Conway}.
Graphs of this type can be enumerated by methods of Tutte \cite{Tutte}. 
\bigbreak

In knot theory, each node of a $4$-regular plane graph is endowed with extra structure so that the graph can be seen to represent
a projection of a curve that is embedded in three dimensional space. We refer to \cite{Kauffman,Harary} for this terminology.
The present paper will be concerned only with the plane graphs themselves, but we shall call a graph that could be the shadow of a knot or link
a {\em knot graph} or a {\em link graph}, respectively. A knot graph is a $4$-regular plane graph such that one can traverse the entire graph by
entering and leaving each node along edges that are {\em not} adjacent in the cyclic order of the local edges at this node. In a link graph there 
will be several such circuits, corresponding to the components of the link's embedding in three-dimensional space.
\bigbreak

We shall use the words {\em vertex, node} and {\em crossing} as synonyms in the discussion that follows.
\bigbreak

Let $G$ be a connected plane graph. Then $G$  intrinsically has $v=v(G)$ nodes,
$e=e(G)$ edges and, from its embedding in the plane, $f=f(G)$ faces.  The faces of
$G$ are the connected regions into which the complement of G in the plane (or in the 
surface of the sphere) is divided. Euler proved the fundamental result:
\bigbreak

\noindent
{\sc Theorem 1.} {\it Let $G$ be a connected plane graph. Then}
$$v(G) - e(G) + f(G) = 2.$$
\bigbreak

\noindent
{\sc Proof.} See \cite{Berge}.  $\hfill \Box$ 
\bigbreak

Let $G$ be a lune-free graph (by which we mean a simple $4$-regular connected
plane graph). Four-regularity implies that $4v(G) = 2e(G)$. Hence,
$$2v(G) = e(G).$$

\noindent Substituting this equation into Euler's formula, we obtain the equation
$$f(G) = 2 + v(G)$$
\noindent when $G$ is a connected lune-free graph. 
\bigbreak

If $G$ is any plane graph, we let $f_{k}(G)$ denote the number of faces of $G$ that 
have $k$ edges on the boundary. Then we can write
$$f(G) = f_{3}(G) + f_{4}(G) + f_{5}(G) + \cdots$$
when $G$ is simple and lune-free, since in this case $f_{1}(G)=f_{2}(G) = 0.$

\bigbreak

\noindent {\sc Theorem 2.} {\it Let $G$ be a connected lune-free link graph. Then}
$$f_{3}(G) = 8 + f_{5}(G) + 2f_{6}(G) + 3f_{7}(G) + \cdots.$$

\noindent {\it Thus any lune-free graph must have at least eight three-sided regions.}
\bigbreak

\noindent {\sc Proof.}
Let $f=f(G), f_k = f_{k}(G), v=v(G), e= e(G).$ We have
$$f = f_3 + f_4 + f_5 + \cdots,$$
$$2e = 3f_3 + 4f_4 + 5f_5 + \cdots,$$
$$ 4f = 2e + 8.$$

\noindent Hence

$$(4f_3 + 4f_4 + 4f_5 + \cdots) = (3f_3 + 4f_4 + 5f_5 + \cdots) + 8,$$

\noindent from which the theorem follows at once. $\hfill \Box$ 
\bigbreak

It follows from Theorem 2
that the smallest (meaning the smallest number of nodes) possible lune-free link graph has $f_3 = 8$ and $f_k = 0$ for all $k$
greater than three. See Figure 1. This diagram is the familiar three-circle Venn diagram,
and it is the unique minimal lune-free link graph. Note that the number of link components
of this diagram is three. We are independently interested in the number of 
lune-free link graphs with a given number of link components.
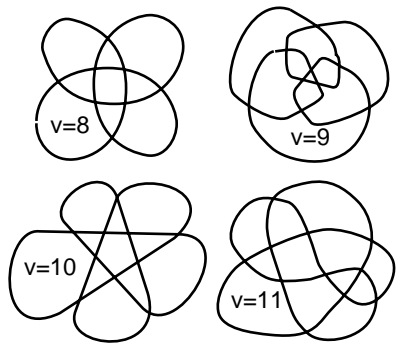
\bigbreak

$$ \picill1inby1.3in(Venn)  $$
\begin{center}
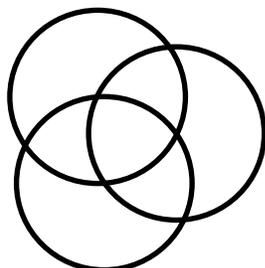
{\sc Figure 1 - The Simplest Lune-Free Link Graph}
\end{center}
\bigbreak

$$ \picill2inby1.8in(LFD911)  $$
\begin{center}
{\sc Figure 2 - Table of Lune-Free Knot Graphs with $8, 9,10$ and $11$ Nodes}
\end{center}
\bigbreak

The next simplest link graph has $f_3=8$ but $f_4 = 2$ and (necessarily)
$f_k = 0$ for all $k > 4.$ To see this, start with a four sided region and begin
building three sided regions adjacent to it. Four more three sided regions are then
forced, and one obtains closure with one extra four sided region on the outside.
This graph is minimal and unique with specifications as above. It is an
eight crossing knot graph with one knot-theoretic component, hence a lune-free 
knot graph with eight vertices. See Figure 2.
\bigbreak

In Figures 2 and 3 we give tables of lune-free knot graphs with
nine through twelve crossings.
\bigbreak

$$ \picill2inby2.2in(LFD12)  $$
\begin{center}
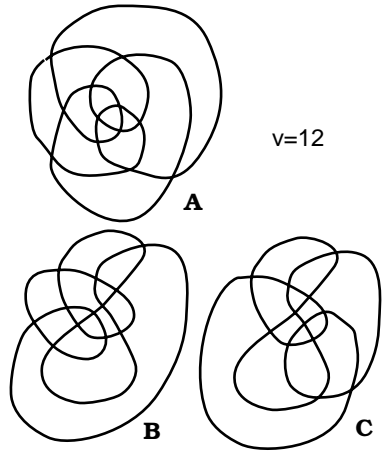
{\sc Figure 3 - Table of Lune-Free Knot Graphs with $12$ Nodes}
\end{center}
\bigbreak

Here is a sketch of the rest of this section and the next section: In the remainder of this section we give various examples including
link graphs with eight triangle faces and all other faces of order four. We begin by giving a bit more information about the Venn diagram and the
graph $G_8$ of Figure 5. There
are infinitely many such examples. In Section 3 we show that there exist lune-free knot graphs
with $v$ any integer greater than or equal to 8.
\bigbreak

\noindent {\sc Example 1.} A smallest lune-free graph must have $f_{3}(G) =8.$ Witness the standard $3-$circle
Venn diagram as shown in Figure 1. This diagram is the unique lune-free graph with $f_{3} = 8,$ as is easy to see
from the construction in Figure 4.

$$ \picill2inby1.2in(Construct)  $$
\begin{center}
{\sc Figure 4 - Creating the Venn Diagram}
\end{center}
\bigbreak

$$ \picill2inby1.2in(G)  $$
\begin{center}
{\sc Figure 5 - The Graph $G_{8}$}
\end{center}
\bigbreak

\noindent {\sc Example 2.} Let $G_{8}$ denote the graph shown in Figure 5.

\noindent Here we have $f_{3} = 8, f_{4} =2, f = 10, v = 8, e = 16,$ and $G_{8}$ is a knot projection (the 
projection of a single component link). This graph is the next largest (in size the number of nodes $v(G)$) lune-free 
graph after the Venn diagram of Example 1. To see this, note that any lune-free graph with an $n$-sided region must have
at least $2n$ nodes. Thus a smallest possible lune-free graph with $f > f_{3}$ must have as least $8$ nodes. $G_{8}$ is the 
smallest such example. Note also that the Venn diagram is the unique lune-free graph with $f=f_{3}.$
\bigbreak

$$ \picill3inby3in(GFamily)  $$
\begin{center}
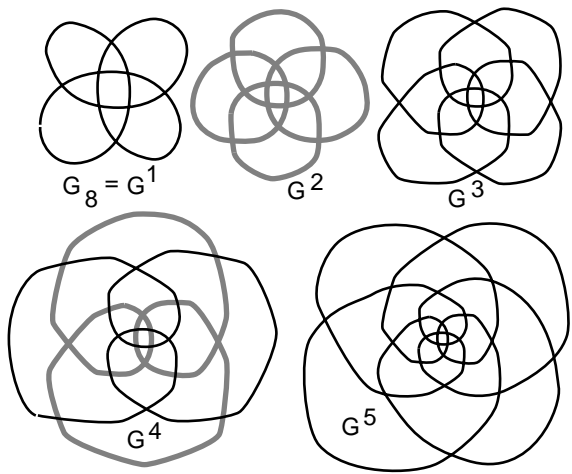
{\sc Figure 6 - The Graph $G_{8}$ and Its Family}
\end{center}
\bigbreak

\noindent {\sc Example 3.} $G_{8}$ is the smallest member of an infinite family of lune-free graphs, as illustrated in
Figure 6.

\noindent In this family, the graphs $G^{n}$ have a single knot theoretic component for $n$ odd. Thus in the Figure
$G^{1} = G_{8}, G^{3}, G^{5}$ are knot graphs. The graph $G^{n}$ has $4n + 4$ nodes.
\bigbreak

\noindent {\sc Example 4.} We can regard the family $\{ G^{n} \}$ as generated from a $4$-gon (the central square in each
figure). A similar family is generated from any $n$-gon. For example, if we start with $n=5,$ the smallest example is the 
graph shown in Figure 7. Note that this pentagonal example has ten nodes.

$$ \picill3inby2in(PentagonFamily)  $$
\begin{center}
{\sc Figure 7 - The Pentagon Family}
\end{center}
\bigbreak

\noindent In the same vein, we can start with a triangle and look at the family that starts with the Venn diagram.
This family is illustrated in Figure 8. In this case, we obtain knot graphs for $9$ nodes, $12$ nodes and $18$ nodes.
For this series
we have the number of nodes $v = 3n$ for $n = 1,2,3,4,...$ where the Venn diagram is the case $n=2.$ (In this mode of counting, $n=1$  is the
trefoil knot graph. Since the trefoil has lunes, we start this series at $n=2.$) Let $\mu(G)$ denote the number of knot theoretic components of a 
knot graph $G.$ Let $\mu(n)$ denote the number of knot theoretic 
components of the $n$-th graph in the trefoil series. The component count for this series is given by the formula
$\mu(n) = 3$ if $n = 2 + 3k$ and $\mu(n) = 1$ if $n = 3k$ or $n = 3k +1.$  The result is easily verified by induction.
\bigbreak
 
$$ \picill3inby3in(VennFamily)  $$
\begin{center}
{\sc Figure 8 - The Venn Family}
\end{center}
\bigbreak

\noindent {\sc Example 5.} At this point we have constructed lune-free knot graphs ($\mu = 1$) for 
nodes $v = 8,9,10,12$ (and some higher order nodes). An $11$ node knot graph appears as the modification of
the graph for $v = 10,$ as shown in Figure 9. This device produces new lune-free graphs from 
lune-free graphs, but does not always preserve the knot component count $\mu.$
\bigbreak

$$ \picill3inby2in(CrossingDevice)  $$
\begin{center}
{\sc Figure 9 - Crossing Device}
\end{center}
\bigbreak

\noindent {\sc Example 6.} At this stage we are prepared to exhibit, in Figures 2 and 3, a complete list of
all $\mu = 1$ lune-free graphs $G$ with less than or equal to twelve nodes. It is interesting to note that the 
graph of 9 nodes in Figure 2 can be obtained from the graph of 8 nodes in Figure 2, by the method of Figure 9.
There are three distinct lune-free graphs with $\mu = 1$ and $v = 12$ nodes. Note that in Figure 3,
the three distinct lune-free graphs with 12
nodes are  labeled $A$, $B$ and $C.$  Graphs $A$ and $B$ each have only 3 and 4 sided regions, but are not isomorphic graphs.
Graph $C$ has a five-sided region, making it distinct from the other two.  
\bigbreak
In Figure 10, we show the complete list of lune-free graphs with 13 nodes. We thank Slavik Jablan for the use of his tables
of Conway polyhedra \cite{Jablan} for this information.
\bigbreak

$$ \picill3inby3.3in(LFD13)  $$
\begin{center}
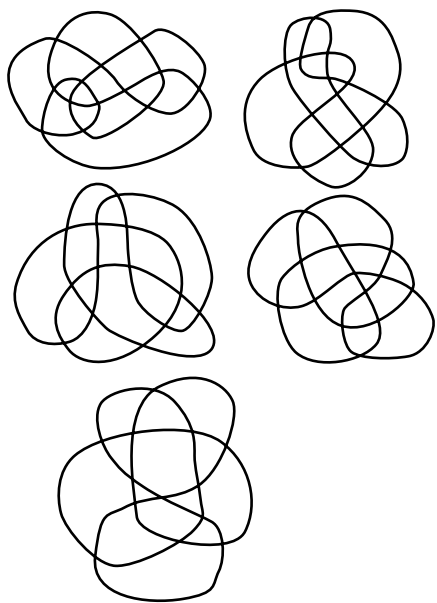
{\sc Figure 10 - All Lune-Free Graphs with Thirteen Nodes}
\end{center}
\bigbreak

\section {Lune-Free Graphs of Every Order}
We now turn to the following question: Is there a lune-free knot graph for each number of
nodes $v \ge 8?$ We shall prove that the answer to this question is {\it yes} by an inductive construction.
\smallbreak

\noindent {\sc Theorem 3.} {\it Let $v \ge 8.$ Then there exists at least one lune-free knot graph with
$v$ nodes.}
\smallbreak

\noindent {\sc Proof.} First note the inductive construction indicated in Figure 11.
In Figure 11, the letters $a,b,c,d,e,f$ denote the number of edges in the boundary of the regions in which
they reside. The graph $G'$ is obtained from the graph $G$ as illustrated. In $G'$ the $a,b,e,f$ regions
are incremented to $a+1, b+1, e+1, f+1$ respectively, while the $c$ and $d$ regions are decremented to 
$c-1$ and $d-1.$ Two new three-sided regions are produced in the process.  Thus the graph $G'$ is lune-free
exactly when $c \ge 4$ and $d \ge 4.$  In other words, the process $G \longrightarrow G'$ induces a lune-free $G'$
whenever $G$ is lune-free and has two adjacent regions with each region having at least $4$ sides. When $G$ satisfies this condition, we shall
call $G$ {\it admissible}. It follows from this construction that if $G$ is admissible, then $G'$ is admissible.
If $v(G) = v,$ then $v(G') = v + 2.$ Hence, for any lune-free admissible $G$ with $v(G)=v,$ we obtain an infinite
family of lune-free graphs $\{ G', G'', G''', ... , G^{(n)}, ... \}$ with $v(G^{(n)}) = v + 2n.$
\smallbreak

Now examine the library of lune-free graphs in Figures 2 and 3. The graphs for $v = 8,9,10$ are not admissible. The
graph for $v=11$ is admissible, and each of the graphs for $v=12$ is admissible. Hence the two inductive processes
starting from $11$ nodes and from $12$ nodes together produce lune-free graphs of all orders greater than 12. Since we have 
already given examples from order $8$ to order $12,$ this completes the proof of the theorem. $\Box$
\bigbreak

$$ \picill3inby2in(IC)  $$
\begin{center}
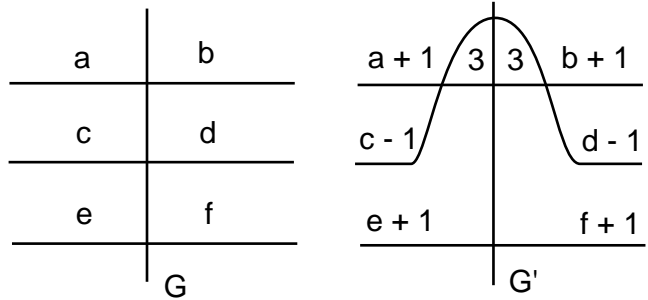
{\sc Figure 11 - Inductive Construction for Increasing Nodes by Two}
\end{center}
\bigbreak

$$ \picill3inby2in(Process)  $$
\begin{center}
{\sc Figure 12 - Illustration of Inductive Process}
\end{center}
\bigbreak

\noindent {\sc Remark.} A series of lune-free graphs of orders $11,13,15,17, \cdots$ obtained by the $G \longrightarrow G'$
process is illustrated in Figure 12.
\bigbreak

\noindent {\sc Remark.} Other inductive constructions are available that can produce lune-free knot projections from given ones.
Some constructions do not depend upon any admissibility condition. One example of such a construction is shown in Figure 13.
In this construction, the new graph has order $v + 9$ where $v$ is the number of nodes in the original graph.
\bigbreak

$$ \picill3inby3in(SC)  $$
\begin{center}
{\sc Figure 13 - From $v$ to $v + 9.$}
\end{center}
\bigbreak

\noindent {\sc Remark.} 
The referee reported the following suggestion, from a colleague,
of a possibly simpler proof of Theorem 3: just consider the shadows of some
explicit 3-braid closures, namely
$$
(\sigma_1 \sigma_2 \sigma_3)^{4k-1} (\sigma_2 \sigma_1)^m \sigma_2^l,
$$
with $k,m,l$ integers satisfying $k \ge 1$, $0 \leq m \leq 5$ and $0
\leq l \leq 1$.

This is an interesting attempt, but it doesn't fully work as is. Indeed,
denote by $S(k,m,l)$ the corresponding link graph. A careful
analysis reveals that $S(k,m,l)$ is not a \textit{knot} graph if $m=2$,
$m=5$ or $m \equiv l \equiv 1 \bmod 3$. In summary, the proposed
construction only gives lune-free knot graphs on $v = 3(4k-1) +2m+l$
nodes if $v \equiv 3,4,5,9,10$ or $11 \bmod 12$.
Moreover, the graph $S(k,m,l)$ is tight (see below) only if $m = l = 0$,
i.e. if $v=12k-3$.
\bigbreak

\noindent {\sc Definition.} A $4$-regular plane graph is said to be {\it tight} if for every pair of adjacent faces, one face
has $3$ edges. Thus tight is synonymous with {\it inadmissible} in the terminology above.
In the next section, we shall consider questions about tight lune-free knot projections.
\bigbreak

\section {Tight lune-free knot graphs}

In the preceding section, we have shown that there is an {\em admissible}
lune-free knot graph with $v$ crossings for every $v\geq 11.$
Here we address the existence question of {\em tight}, i.e.
inadmissible, lune-free knot and link graphs.
In the above tables, we have already encountered such graphs
for $v= 8,9$ and 10. The basic question, that we shall completely
answer here, is the following: for which integers $v$ does there
exist a {\em tight} lune-free knot or link graph, respectively?
For links, the answer is: all integers $v \geq 8$ such that
$v \not \equiv \pm 1 \bmod 6$ can be realized as the number of crossings
of such graphs. The answer for knots is the same as for links, with the sole 
exception of $v=12$, which cannot be realized by a tight lune-free {\em knot} 
graph. However, here is a tight lune-free link diagram with 12 crossings and 
$\mu = 3$ components.
\bigbreak

$$ \picill3inby1.8in(Page1)  $$
\begin{center}
{\sc Figure 14 - A Symmetrical Example with $v=12.$}
\end{center}
\bigbreak

The aim of this section is to prove the following result.

\bigskip
{\sc Theorem 4.} 
{\em (1) There exists a tight lune-free link graph with $v$
crossings if and only if $v \geq 8$ and $v\not \equiv \pm 1 \bmod 6$.
(2) There exists a tight lune-free knot graph with $v$
crossings if and only if $v \geq 8$, $v\not \equiv \pm 1 \bmod 6$, 
and $v \not = 12$.}
\bigskip


For the proof, which will occupy the rest of the section, we shall need 
several preliminaries.
We start by recalling a few things about the {\em medial graph} $D = m(G)$ 
of a plane graph $G$. The vertices of $D$ correspond to the edges of $G$, 
and an edge in $D$ corresponds to a pair of edges in $G$ which lie on a 
common face. Thus, $D= m(G)$ is a 4-regular plane graph, i.e. a link graph.
Figure 15 gives an example.
Conversely, every link graph $D$ arises as the medial graph of a suitable
plane graph $G$, obtained as follows. Color the faces of $D$ in black and
white in checkerboard manner. The vertices of $G$ correspond to the black 
faces of $D$, and there is an edge between two vertices in $G$ for each 
common vertex of the corresponding two black faces in $D$. 
\medskip

Note that the black faces of $D$ correspond to vertices of $G$, while 
the white faces of $D$ correspond to faces of $G$. Thus, 
$m(G) = m(G^{\ast})$, where $G^{\ast}$ denotes the dual of $G$.
From this description, it follows that a pair of adjacent faces in 
$D = m(G)$ corresponds to a pair formed by one vertex and one face 
incident to it in $G$. In other words, the faces of $D=m(G)$ are in
one-to-one correspondence with the vertices and faces of the original
plane graph $G$. See Figure 16.

\bigbreak

$$ \picill3inby2in(Page20)  $$
\begin{center}
{\sc Figure 15 - Medial Graph}
\end{center}
\bigbreak

$$ \picill3inby2in(Page2)  $$
\begin{center}
{\sc Figure 16 - More on the Medial Graph}
\end{center}
\bigbreak

For our purposes here, we do need to address the following two questions. 
Let $G$ be a simple connected plane graph, and let $D=m(G)$ be its 
medial graph. Under what conditions on $G$ can one guarantee that:
\medskip

(1) $D$ is a {\em knot} graph, rather than a many-component link graph?
\medskip

(2) $D$ is a {\em tight} lune-free graph?
\bigskip

In order to answer question (1), we shall make the following definitions.
Let $G$ be a simple connected plane graph.
\bigskip

{\sc Definition.} 
An {\em angle} in $G$ is a pair $\alpha = (v,F)$ where
$v$ is a vertex of $G$ and $F$ is a face of $G$ containing $v$.
See Figure 17.
\bigbreak

$$ \picill3inby1in(Page31)  $$
\begin{center}
{\sc Figure 17 - The angle $(v,F)$}
\end{center}
\bigbreak

$$ \picill3inby1.5in(Page32)  $$
\begin{center}
{\sc Figure 18 - Two adjacent angles}
\end{center}
\bigbreak

{\sc Definition.} Two angles $\alpha_1 = (v_1,F_1)$ and $\alpha_2 = (v_2,F_2)$
are {\em adjacent} if $v_1,v_2$ are connected by an edge $e$, with
$e = F_1 \cap F_2$. Figure 18 is
an illustration of adjacency.
\smallbreak

One may thus form  the {\em graph ${\cal A}(G)$ of angles of
$G$.} Clearly, each angle is adjacent to exactly two
other angles in $G$. It follows that the graph ${\cal A}(G)$
is regular of degree 2, {\it i.e.} ${\cal A}(G)$ is a disjoint
union of circuits. We may now give an answer to question (1)
above.
\bigskip

{\sc Proposition 5.} {\em Let $G$ be a simple connected plane graph,
and let $D=m(G)$ be its medial graph. Then the number of
components of any link represented by the link graph $D$
is equal to the number of connected components of the
graph of angles ${\cal A}(G)$.}
\bigskip

The proof of this proposition is easy, and will be left to 
the reader. It suffices to think carefully about the way
the medial graph is formed, which essentially amounts to connect
adjacent angles by strands. Note that this result was 
independently observed by Isidoro Gitler as well as Sostenes Lins.
\bigbreak

The referee points out that the Tutte polynomial gives another way to  compute
the number of components $c$ of a link represented by the medial graph of $G$.
However, for practical purposes, Proposition 5 is a very efficient way to compute $c$
and actually guided us for the constructions in our proof below of Theorem 4.
\bigbreak

Figure 19 illustrates an example of a graph $G$ whose medial graph $D=m(G)$
represents 3-component links. In the Figure, the 3 different
sequences of adjacent angles are denoted by 
$a,b,c$, respectively.
\bigbreak

$$ \picill3inby1.5in(Page4)  $$
\begin{center}
{\sc Figure 19 - Three paths of adjacent angles }
\end{center}
\bigbreak

We shall now address question (2). The following definitions
will be helpful for this purpose.
\bigskip

{\sc Definition.} Let $F$ be a face in a given plane graph $G$. The 
{\em degree} of $F$, denoted deg($F$), is the number of edges of $G$ 
which are incident to $F$.
\bigbreak
\bigbreak

{\sc Definition.} A connected plane graph $G$ is {\em special} if 
\begin{enumerate}
\item all vertices and faces of $G$ have degree $\geq 3$
\item if a vertex $x$ is incident with a face $F$, then either
deg($x$)=3 or deg($F$)=3. 
\end{enumerate}
\bigskip

{\sc Lemma 6.} {\em Let $G$ be a simple connected plane graph,
and let $D=m(G)$ be its medial graph. Then $D$
is a tight lune-free link graph if and only if $G$
is special in the above sense.}
\bigskip

{\em Proof.} Clearly, the lune-free condition for $D$ corresponds
to condition (1) on $G$ as a special graph. Similarly,
condition (2) on $G$ exactly means that $D=m(G)$ cannot have 
two adjacent faces of degree $\geq 4$ each. $\Box$
\bigskip

Examples of special graphs include the connected plane cubic graphs,
the sphere triangulations (which are the duals of plane cubic graphs), 
and the wheels. All these graphs are easily seen to have a medial graph 
which is a tight lune-free link graph. We shall now show that there 
are no other special graphs.
\bigskip

{\sc Theorem 7.} {\em Let $G$ be a simple connected plane graph. Then $G$ 
is special if and only if $G$ is either cubic, or a triangulation, or a 
wheel.}
\bigskip

{\em Proof.} We have already observed that cubic graphs, sphere
triangulations and wheels are special graphs. Conversely,
assume that $G$ is special, and that it is neither cubic
nor a triangulation. We shall prove that $G$ is then necessarily
a wheel.
\bigbreak

Let $L$ (for large) be the set of vertices of degree $\geq 4$, 
and $S$ the set of vertices which are adjacent to a face
of degree $\geq 4$. Neither $L$ nor $S$ is empty, as $G$
is neither cubic nor a triangulation. Moreover, since $G$ is special, 
all vertices in $S$ have degree exactly 3, i.e. $L \cap S = \emptyset$.   
Let $d \geq 1$ be the minimal distance between vertices in $L$ and 
vertices in $S$, and let $x_0 \in L$, $y_0 \in S$ at distance exactly 
$d$.
Let $W$ be the subgraph induced by $x_0$ and its neighbors. 
Then $W$ is a wheel with center $x_0$,  because all faces adjacent 
to $x_0$ must have degree 3. We will show that $G$ = $W$, thereby
completing the proof of the theorem.
\bigbreak

Let $x_0, x_1, \ldots, x_{d-1}, x_d = y_0$ be a path in $G$ of length $d$
joining $x_0$ to $y_0$. Of course, $x_1$ belongs to $W$. We must have
deg($x_1$) = 3, for otherwise $x_1$ would belong to $L$, contradicting
the minimality of $d$. Note that all three neighbors of $x_1$ lie in $W$.
We will call $z_1$, $z_2$ the two neighbors of $x_1$ which are distinct 
from $x_0$. 
Assume for a contradiction that $d\geq 2$. Then $x_2$ does not belong
to $W$, for otherwise there would be a shorter path from $x_0$ to $y_0$.
This contradicts the fact that all neighbors of $x_1$ belong to $W$.
Therefore, we must have $d=1$, i.e. $x_1=y_0$. Let $F$ be the face 
of degree $\geq 4$ to which $y_0$ is adjacent, and let $C$ be the cycle
bounding $F$. Since $G$ is special, all vertices of $C$ have degree exactly 3.
Now of course $x_1$, $z_1$ and $z_2$ belong to $C$. Hence all neighbors
of $z_1$ or $z_2$ already lie in $W$. As easily seen, this implies that 
$C$ is entirely contained in $W$, and actually coincides with the boundary 
of $W$. 
Thus, all vertices of $W$ except its center $x_0$ have degree exactly three. 
Therefore $G$ = $W$, for otherwise, since $G$ is connected, there 
would be vertices in the boundary of $W$ of degree higher than 3. 
It follows that $G$ is a wheel, as claimed. Figure 20 illustrates 
the path combinatorics used in this proof. $\hfill \Box$

\bigbreak

$$ \picill3inby1in(Page6)  $$
\begin{center}
{\sc Figure 20 - Connecting Path}
\end{center}
\bigbreak

{\sc Corollary 8.} {\em Let $D$ be a tight lune-free link graph. Let
$G$ be a simple connected plane graph whose medial graph $m(G)$ is
equal to $D$. Then $G$ is either a cubic graph, a sphere triangulation,
or a wheel.}  $\Box$
\bigskip

\bigskip

We shall now turn to the proof of the main theorem.
\bigskip

{\em Proof of Theorem 4.} 
\medskip

(a) {\em The case $v \equiv \pm 1 \bmod 6$.}
Our first step will be to prove that there cannot be a tight
lune-free link graph with $v$ crossings if 
$v \equiv \pm 1 \bmod 6$.
Assume the contrary, and let $D$ be such a graph. Let $G$
be a plane graph such that $m(G) = D.$ From the hypotheses on
$D$ and the above lemma, it follows that $G$ is special graph.
Now, by the above theorem, $G$ must be either cubic, or a triangulation, 
or a wheel. On the one hand, if $G$ is cubic or a triangulation, 
then the number of edges in $G$ is a multiple of 3. Thus 
$v \equiv 0$ or $3 \bmod 6$ in these two cases. On the other hand, 
if $G$ is a wheel, then it has an even number of edges, whence 
$v \equiv 0, 2$ or $4 \bmod 6$ in this third and last case. 
\medskip

(b) {\em The other cases.}
We shall now construct tight lune-free {\em knot} graphs 
with $v$ crossings, whenever $v \geq 8$, $v \equiv 0,2,3$ or 
$4 \bmod 6$, and $v \not = 12.$
This will be done by induction, and the basic induction step 
is provided by the following result.
\bigskip

{\sc Proposition 9.} {\em If there is a tight knot graph with $v$ 
crossings, then there is a tight knot graph with $v+12$ crossings.}
\bigskip

{\em Proof.} Let $G$ be a plane graph with $v$ edges, such that 
its medial graph $D=m(G)$ is a tight knot graph with $v$
crossings. In particular, $G$ is a special graph.
Since $D$ has at least one (in fact at least 8) 3-face, we may assume,
taking the dual of $G$ if necessary, that $G$ has a vertex of degree 3.
Now do the move shown in Figure 21 on $G$ near this vertex, calling $H$ the 
resulting graph.
\bigbreak

$$ \picill3inby1.5in(Page70)  $$
\begin{center}
{\sc Figure 21 - Laddering }
\end{center}
\bigbreak

Clearly, $H$ has $v+12$ edges, and it is obvious that $H$ is special
since $G$ is. Thus, $m(H)$ is a tight {\em link} graph with $v+12$
vertices. 

It remains to show that $m(H)$ is in fact a {\em knot} graph. To do
this, it suffices to compare the paths of adjacent angles in $G$ 
and in $H$, and observe that the configurations are completely similar.
(See Figures 22 and 23.) Since $G$ yields a knot, so will $H$.

In $G$, we picture three partial paths of adjacent angles,
namely $(a_1,a_2,a_3)$, $(b_1,b_2)$ and $(c_2,c_3)$. (See Figure 22.) 

$$ \picill3inby1.5in(Page7)  $$
\begin{center}
{\sc Figure 22 - Adjacent angles in $G$}
\end{center}
\bigbreak

$$ \picill3inby3.2in(Page71)  $$
\begin{center}
{\sc Figure 23 - Corresponding adjacent angles in $H$}
\end{center}
\bigbreak

In $H$, the
corresponding partial paths are $(a_1,a_2,\dots,a_{15})$, 
$(b_1,\dots,b_8)$ and $(c_2,\dots,c_9)$. (See Figure 23.) Thus, $m(H)$ is a 
knot graph, as asserted. $\Box$
\bigbreak

$$ \picill3inby1.7in(Page8)  $$
\begin{center}
{\sc Figure 24 - Weaving a Medial Graph}
\end{center}
\bigbreak
\bigbreak

We shall now complete the proof of Theorem 4.
\medskip

(1) The case $v \equiv 2,4$ mod 6.
Let $n = v/2$ and $G = W_n$ be the wheel with $n+1$ vertices. 
Then $m(G)$ is a tight lune-free link graph with $v=2n$ vertices, 
and is in fact a knot graph as $n$ is not divisible by 3. 
\smallbreak

(2) The case $v \equiv 0,3$ mod 6.
Using the $v \mapsto v + 12$ inductive construction in the above
proposition, it suffices to exhibit tight lune-free knot graphs 
with $v= 9,15,18$ and $24$ crossings respectively.
These are shown in Figure 25. $\Box$
\smallbreak

$$ \picill3inby2.3in(Page9)  $$
\begin{center}
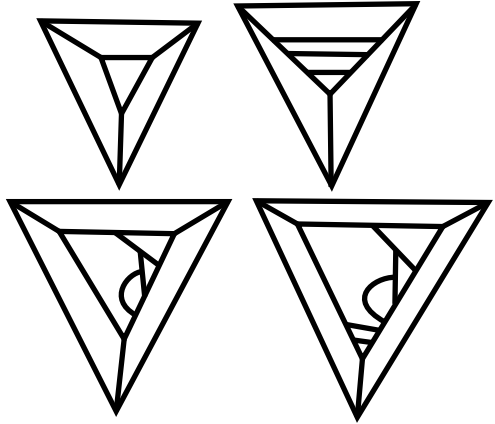
{\sc Figure 25 - Plane graphs whose medials are tight lune-free knot graphs with $v= 9,15,18$ and $24$ }
\end{center}

\bigbreak

\section{Knot graphs with exactly $k$ lunes.}
\bigbreak

Let $k \geq 0$ be an integer. A knot graph is {\bf $k$-lune} if it has
exactly $k$ lunes. We will show that if $v$ is large enough with respect
to $k$, then there is a $k$-lune knot graph with $v$ crossings.
The case $k=0$ is equivalent to being lune-free.
\bigbreak

{\sc Proposition 10.} {\em Let $k \geq 0$ be an integer. If $k$ is even
and $v \geq k+8$, or if $k$ is odd and $v \geq k+9$, then there is
a $k$-lune knot graph with $v$ crossings.}
\bigbreak

{\em Proof.} 
The case $k = 0$, i.e. the lune-free case, has been dealt with in an 
earlier section.
\bigbreak

(0) Let $k$ even, $k \geq 2$. Let $v \geq k+8$, and let $D$
be a lune-free knot graph with $v-k$ crossings, which exists
as we already know.
Necessarily, $D$ has a 3-face (in fact at least 8, by Euler's formula
and Theorem 2). Now perform the 
deformation of this 3-face as in Figure 26.
\bigbreak

$$ \picill3inby1.2in(Page10)  $$
\begin{center}
{\sc Figure 26 - 3-face deformation}
\end{center}
\bigbreak

\noindent
Add $k$ new crossings to $D$, with $k/2$ lunes above and $k/2$
lunes below the desired edge. The new knot graph has exactly $k$
lunes and $v$ crossings.
\medskip

(1) The case $k$ = 1. Let $v \geq 10$ and let $D$ be a lune-free knot 
graph with $v-2$ crossings. Then $D$ admits a large face, say $F$.
Take two non-adjacent edges bounding $F$, and do the move shown in Figure 27.
\bigbreak

$$ \picill3inby1in(Page101)  $$
\begin{center}
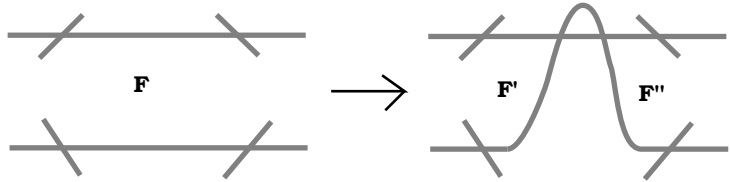
{\sc Figure 27 - Move on Two Non-Adjacent Edges}
\end{center}
\bigbreak

The new knot graph on the right has exactly 1 lune and $v$ crossings,
as desired.
\medskip

(2) Let finally $k$ be odd, with $k \geq 3$. Let $v \geq k+9$, and let
$D$ be a knot graph with $v-k+1$ crossings and exactly one lune, 
which exists by the above case (1) and the fact that $v-k+1 \geq 10$.
Then $D$ has a 3-face not adjacent to its unique lune. The same 
modification of this 3-face as in case (0) above, adding $k-1$
new crossings, $(k-1)/2$ lunes above and $(k-1)/2$ lunes below
an edge of this 3-face, will produce a knot graph with $k$
lunes and $v$ crossings. $\Box$

 \end{document}